\documentclass[preprint,12pt]{elsarticle}

\usepackage{amssymb}
\usepackage{amsmath}
\usepackage{amsthm}
\usepackage[colorlinks]{hyperref}
\AtBeginDocument{%
  \hypersetup{%
    linkcolor=blue,%
    citecolor=green,%
  }%
}
\usepackage{cleveref}
\usepackage{geometry}
\usepackage{amsthm}
\usepackage{upgreek}

\crefname{equation}{}{}

\newtheorem{theorem}{Theorem}[section]
\newtheorem{lemma}[theorem]{Lemma}

\theoremstyle{definition}
\newtheorem{definition}[theorem]{Definition}

\theoremstyle{remark}
\newtheorem{remark}[theorem]{Remark}

\numberwithin{equation}{section}

\journal{~}

\begin{document}

\begin{frontmatter}



\title{Boundary Lipschitz Regularity and the Hopf Lemma on Reifenberg Domains for Fully Nonlinear Elliptic Equations \tnoteref{t1}}

\author[rvt]{Yuanyuan Lian}
\ead{lianyuanyuan@nwpu.edu.cn; lianyuanyuan.hthk@gmail.com}
\author[rvt]{Wenxiu Xu\corref{cor1}}
\ead{wenxiuxu@mail.nwpu.edu.cn}
\author[rvt]{Kai Zhang}
\ead{zhang\_kai@nwpu.edu.cn; zhangkaizfz@gmail.com}
\tnotetext[t1]{This research is supported by the National Natural Science Foundation of China (Grant No. 11701454) and the Natural Science Basic Research Plan in Shaanxi Province of China (Program No. 2018JQ1039).}

\cortext[cor1]{Corresponding author}

\address[rvt]{Department of Applied Mathematics, Northwestern Polytechnical University, Xi'an, Shaanxi, 710129, PR China}

\begin{abstract}
In this paper, we prove the boundary Lipschitz regularity and the Hopf Lemma by a unified method on Reifenberg domains for fully nonlinear elliptic equations. Precisely, if the domain $\Omega$ satisfies the exterior Reifenberg $C^{1,\mathrm{Dini}}$ condition at $x_0\in \partial \Omega$ (see \Cref{d-re}), the solution is Lipschitz continuous at $x_0$; if $\Omega$ satisfies the interior Reifenberg $C^{1,\mathrm{Dini}}$ condition at $x_0$ (see \Cref{d-H-re}), the Hopf lemma holds at $x_0$. Our paper extends the results under the usual $C^{1,\mathrm{Dini}}$ condition.
\end{abstract}


\end{frontmatter}


\section{Introduction}\label{S1}

In this paper, we intend to obtain the pointwise boundary Lipschitz regularity and prove the Hopf Lemma for the viscosity solutions of the following fully nonlinear elliptic equations
\begin{equation}\label{r-e}
\left\{\begin{aligned}
&u\in S(\lambda,\Lambda,f)&& ~~\mbox{in}~~\Omega;\\
&u=g&& ~~\mbox{on}~~\partial \Omega,
\end{aligned}\right.
\end{equation}
where $\Omega$ is a bounded domain and $S\left( {\lambda ,\Lambda ,f} \right)$ denotes the Pucci class with uniform constants $\lambda$ and $\Lambda$ (see \cite{MR1351007} for the definition and its basic properties).

It is well known that the exterior sphere condition and the interior sphere condition imply the boundary Lipschitz regularity and the Hopf lemma respectively. In recent decades, sphere condtion has been extended to a more general geometrical condition, i.e., the $C^{1,\mathrm{Dini}}$ condition (see \Cref{r-2}). With aid of the boundary Harnack inequality, Safonov \cite{Safonov2008} proved the boundary Lipschitz regularity under the exterior $C^{1,\mathrm{Dini}}$ condition, and the Hopf lemma under the interior $C^{1,\mathrm{Dini}}$ condition for classical solutions of linear elliptic equations in nondivergence form. Huang, Li and Wang \cite{MR3167627} also obtained the boundary Lipschitz regularity for linear elliptic equations under the exterior $C^{1,\mathrm{Dini}}$ condition. They used an auxiliary barrier function and the iteration technique, without using the boundary Harnack inequality. Lieberman \cite{MR779924} proved the Hopf lemma for linear elliptic equations under the interior $C^{1,\mathrm{Dini}}$ condition by applying the regularized distance. Recently, Lian and Zhang \cite{lian2018boundary} extend above results to fully nonlinear elliptic equations by a unified method. Moreover, the proof is simple. In that paper, curved boundaries were regarded as the perturbation of a hyperplane. Then the desired regularity can be obtained by a perturbation argument.

In this paper, we prove the boundary Lipschitz regularity and the Hopf lemma for viscosity solutions of fully nonlinear elliptic equations under the Reifenberg $C^{1,\mathrm{Dini}}$ condition which extends the $C^{1,\mathrm{Dini}}$ condition. We use an improved technique of \cite{lian2018boundary} to derive our results. It can be shown that the directions at different scales converges to a direction, say $e_n$. The main difficulty is that we can not estimate the derivative of the solution along $e_n$ directly, which is carried out in  and the difference between two adjoint scales has to be estimated. The Reifenberg $C^{1,\mathrm{Dini}}$ condition was introduced by Ma, Moreira and Wang \cite{MR1351008}where the boundary $C^{1}$ regularity was obtained for fully nonlinear parabolic equations. Note that the Reifenberg ${C^{1,\mathrm{Dini}}}$ condition is more general than the ${C^{1,\mathrm{Dini}}}$ condition (see \Cref{d-re} and \Cref{d-H-re}).

Before the statement of our main results, we introduce some standard notations and definitions. Let $B_r(x_0)$ denote the open ball in $R^n$ with center $x_0$ and radius $r$. Set $B_r=B_r(0)$, $B_r^+=B_r\cap \left\{x|x_n>0\right\}$ and $T_r=B_r\cap \left\{x|x_n=0\right\}$. Denote by $Q_r(x_0)$ the open cube in $R^n$ with center $x_0$ and side-length $r$. Set $Q_r=Q_r(0)$ and $Q_r^+=Q_r\cap \left\{x|x_n>0\right\}$. In this paper, $\left\{e_1,...,e_n\right\}$ stands for the standard basis in $R^n$.

\begin{definition}\label{d-Dini}
The function $\omega:[0,+\infty)\rightarrow [0,+\infty)$ is called a Dini function if $\omega$ is nondecreasing and satisfies the following Dini condition for some $r_0>0$
\begin{equation}\label{e-dini}
    \int_{0}^{r_0}\frac{\omega(r)}{r}dr<\infty.
\end{equation}
\end{definition}

\begin{definition}\label{d-Dini-f}
Let $\Omega \subset R^{n}$ be a bounded domain and $f$ be a function defined on $\bar{\Omega}$. We say that $f$ is Lipschitz at $x_0\in \bar{\Omega}$ or $f\in C^{0,1}(x_0)$ if there exists a constant $C$ such that
\begin{equation*}
  |f(x)-f(x_0)|\leq C|x-x_0|,~~\forall~x\in \bar{\Omega}.
\end{equation*}
Then, define $[f]_{C^{0,1}(x_0)}=\inf C$ where $C$ satisfies above equation, and $\|f\|_{C^{0,1}(x_0)}=\|f\|_{L^{\infty}(\Omega)}+[f]_{C^{0,1}(x_0)}$.

Similarly, we call that $f$ is $C^{1,\mathrm{Dini}}$ at $x_0$ or $f\in C^{1,\mathrm{Dini}}(x_0)$ if there exist a vector $l$ and a constant $C$ such that
\begin{equation*}
  |f(x)-f(x_0)-l\cdot (x-x_0)|\leq C|x-x_0|\omega(|x-x_0|),~~\forall~x\in \bar{\Omega},
\end{equation*}
where $\omega$ is a Dini function. Then we denote $l$ by $\nabla f(x_0)$. We define $[f]_{C^{1,\mathrm{Dini}}(x_0)}=\inf C$ where $C$ satisfies above equation, and $\|f\|_{C^{1,\mathrm{Dini}}(x_0)}=\|f\|_{L^{\infty}(\Omega)}+|l|+[f]_{C^{1,\mathrm{Dini}}(x_0)}$.
\end{definition}

Now, we give the definitions of the Reifenberg $C^{1,\mathrm{Dini}}$ condition.
\begin{definition}[\textbf{Exterior Reifenberg $C^{1,\mathrm{Dini}}$ condition}]\label{d-re}
 We say that $\Omega$ satisfies the exterior Reifenberg ${C^{1,\mathrm{Dini}}}$ condition at $0 \in \partial \Omega$ if there exist a positive constant $R$ and a Dini function $\omega_{\Omega}(r)$ such that

a) for any $0 < r < R,$ there exists a hyperplane ${\Gamma _r}$ and its unit normal vector $n_{\Gamma _r}$ such that
\begin{equation}\label{e-re}
{B_r} \cap \Omega  \subset {B_r} \cap \left\{ {{x \cdot {n_{{\Gamma _r}}}} >  - r\omega_{\Omega}(r)} \right\}.
\end{equation}

b) $\left| n_{\Gamma _r}-n_{\Gamma _{\theta r}} \right| \leq {K(\theta) }\omega_{\Omega}(r)$ for each $ 0<\theta <1$ and $0<r<R$, where the nonnegative function $K$ depends only on $\theta$ and is bounded on $[\theta _0,1]$ for any $0<{\theta _0}<1$.
\end{definition}

\begin{definition}[\textbf{Interior Reifenberg $C^{1,\mathrm{Dini}}$ condition}]\label{d-H-re}
 We say that $\Omega$ satisfies the interior Reifenberg ${C^{1,\mathrm{Dini}}}$ condition at $0 \in \partial \Omega$ if there exist a positive constant $R$ and a Dini functions $\omega_{\Omega}(r)$ such that

a) for any $0 < r < R,$ there exists a hyperplane ${\Gamma _r}$ and its unit normal vector $n_{\Gamma _r}$ such that
\begin{equation}\label{e-re-2}
{B_r} \cap \Omega^c\subset {B_r} \cap \left\{ {{x \cdot {n_{{\Gamma _r}}}} < r\omega_{\Omega}(r)} \right\} .
\end{equation}

b) $\left| n_{\Gamma _r}-n_{\Gamma _{\theta r}} \right| \leq {K(\theta) }\omega_{\Omega}(r)$ for each $ 0<\theta <1$ and $0<r<R$.
\end{definition}

\begin{remark}\label{r-1}
If $\Omega$ satisfies both the exterior and interior Reifenberg ${C^{1,\mathrm{Dini}}}$ condition at $0$, we call that $\partial\Omega$ is Reifenberg ${C^{1,\mathrm{Dini}}}$ at $0$. Without loss of generality, we always assume that $K\geq 1$.
\end{remark}

\begin{remark}\label{r-2}
If $\Gamma_r$ and $n_{\Gamma_r}$ are the same for different $r$, we arrive at the definition of the usual ${C^{1,\mathrm{Dini}}}$ conditions (see \cite[Definition 1.2 and Definition 1.3]{lian2018boundary}). The Reifenberg ${C^{1,\mathrm{Dini}}}$ condition is more general than the usual ${C^{1,\mathrm{Dini}}}$ condition, which is shown clearly by the following example adopted from \cite{MR1351008}.

 Let $\Omega  = {B_1} \cap \left\{ {y > f\left( x \right)} \right\}\subset R^2$ where $f\left( x \right) ={x}/ {\ln \left| x \right|},x \in \left(-1/2,1/2 \right)$. Clearly, $\partial\Omega$ is not ${C^{1,\mathrm{Dini}}}$ at $0$ since $1/|\ln r|$ is not a Dini function. Now we show that it is Reifenberg ${C^{1,\mathrm{Dini}}}$ at $0$.

For any $0<r<1/2$, let $\Gamma_r$ be the line $\left\{(x,y)| y=x/\ln r\right\}$ and take $n_{\Gamma_r}=(-1/\ln r, 1)$. It is easy to see that there exists a unique $-r<x^*<0$ satisfying $f'\left( {{x^*}} \right) = {1}/{\ln r}$, and for any $(x,y)\in B_r\cap \Omega$,

\begin{equation*}
  (x,y)\cdot n_{\Gamma_r}> (x^*,f(x^*))\cdot n_{\Gamma_r}=  \frac{x^*}{\ln |x^*|} - \frac{x^*}{\ln r}.
\end{equation*}
Note that
\begin{equation*}
  f'(x^*)=\frac{1}{\ln |x^*|}-\frac{1}{\ln^2 |x^*|}=\frac{1}{\ln r}.
\end{equation*}
I.e.,
\begin{equation*}
\frac{x^*}{\ln r}=\frac{x^*}{\ln |x^*|}-\frac{x^*}{\ln^2 |x^*|}.
\end{equation*}
Hence,
\begin{equation*}
 (x,y)\cdot n_{\Gamma_r}> \frac{x^*}{\ln^2 |x^*|}\geq \frac{-r}{\ln^2r}.
\end{equation*}
In addition,
\begin{equation*}
|n_{\Gamma_r}-n_{\Gamma_{\theta r}}| \leq \left| {\frac{1}{{\ln r}} - \frac{1}{{\ln \left( {\theta r} \right)}}} \right| \leq \left| {\frac{{\ln \theta }}{{\ln r\ln \left( {\theta r} \right)}}} \right| \leq \left| \frac{\ln \theta }{ \ln^2 r} \right|.
\end{equation*}

Since $1/\ln^2r$ is a Dini function, $\partial\Omega$ satisfies the exterior Reifenberg ${C^{1,\mathrm{Dini}}}$ condition at $0$. Similarly, it can be verified that $\partial\Omega$ also satisfies the interior Reifenberg ${C^{1,\mathrm{Dini}}}$ condition at $0$.
\end{remark}

Now, we state our main results. For the boundary Lipschitz regularity, we have
\begin{theorem}\label{t-2}
Suppose that $\Omega$ satisfies the exterior Reifenberg $C^{1,\mathrm{Dini}}$ condition at $0\in \partial \Omega$ for some Dini function $\omega_{\Omega}(r)$ and $R>0$. Let $u$ be a viscosity solution of
\begin{equation}
\left\{\begin{aligned}
&u\in S(\lambda,\Lambda,f)&& ~~\mbox{in}~~\Omega;\\
&u=g&& ~~\mbox{on}~~\partial \Omega,
\end{aligned}\right.
\end{equation}
where $g$ is $C^{1,\mathrm{Dini}}$ at $0$ with a Dini function $\omega_g$ and $f\in L^{n}(\Omega)$ satisfies
\begin{equation}\label{e-dini-f}
  \int_{0}^{R}\frac{\omega_f(r)}{r}dr:=   \int_{0}^{R}\frac{\|f\|_{L^{n}(\Omega\cap B_r)}}{r\|f\|_{L^{n}(\Omega)}}dr<\infty.
\end{equation}

Then $u$ is $C^{0,1}$ at $0$ and
\begin{equation*}
  |u(x)-u(0)|\leq C |x|\left(\|u\|_{L^{\infty }(\Omega)}+\|f\|_{L^{n}(\Omega)}+[g]_{C^{1,\mathrm{Dini}}(0)}\right), ~~\forall ~x\in \Omega\cap B_{R},
\end{equation*}
where $C$ depends only on $n, \lambda, \Lambda, K,\omega_{\Omega},\omega_{f},\omega_{g}$ and $R$.
\end{theorem}

For the Hopf lemma, we have
\begin{theorem}[\textbf{Hopf lemma}]\label{t-H-2}
Suppose that $\Omega$ satisfies the interior Reifenberg $C^{1,\mathrm{Dini}}$ condition at $0\in \partial \Omega$ for some Dini function $\omega_{\Omega}(r)$ and $R>0$. Let $u\in C(\bar{\Omega})$ satisfy
\begin{equation}
M^{-}(D^2u,\lambda,\Lambda)\leq 0 ~~\mbox{in}~~\Omega~~ (\mathrm{i.e.,}~~ u\in \bar{S}(\lambda,\Lambda,0))
\end{equation}
with $u(0)=0$ and $u\geq 0$ in $\Omega$.

Then there exists a unit vector, say $e_n$, such that for any $l\in R^n$ with $|l|=1$ and $l\cdot {e_n}>0$,
\begin{equation}\label{e-H-main}
  u(tl)\geq cu\left(\frac{e_n}{2}\cdot R \right)t, ~~\forall~ 0<t<\delta,
\end{equation}
where $c>0$ and $\delta>0$ depend only on $n, \lambda, \Lambda,K,\omega_{\Omega},R$ and $l$.
\end{theorem}

\section{Proofs of the main results}
In this section, we give the detailed proofs of the main results.

Now, we clarify the idea briefly. Firstly, note that if the $\Omega$ satisfies the Reifenberg ${C^{1,\mathrm{Dini}}}$ condition from the exterior (or the interior) at $0\in \partial\Omega$, the normal vectors in different scales converges. In addition, the difference of unit normal vectors in different scales is controlled by the Dini function. Next, we use solutions with flat boundaries (i.e., $v$ in the proofs) to approximate the solution $u$. Then the error between $u$ and $v$ (i.e., $w$ in the proofs) can be estimated by maximum principles. By an iteration argument, the boundary regularity for $u$ is obtained. For the boundary Lipschitz regularity, the right hand function $f$, the boundary value $g$ and the curved boundary $\partial \Omega$ are regarded perturbations of $0$, $0$ and a hyperplane (see the definition of $v$ in the proof) which are inspired directly by \cite{MR3780142}. For the Hopf lemma, since the solution is nonnegative and the equation has the right hand zero, it is easier to prove.

First, we prove the following simple lemma.
\begin{lemma}\label{t-1}
Suppose that $\left| n_{\Gamma _r}-n_{\Gamma _{\theta r}} \right| \leq {K(\theta) }\omega_{\Omega}(r)$  for each $0 < \theta  < 1$ and $0<r<R$. Then there exists a unique unit vector $n_0$ such that
\begin{equation}
\mathop {\lim }\limits_{r \to 0 } {{n_{{\Gamma _r}}}} = {n_0} .
\end{equation}
\end{lemma}
\proof For any $0<\eta<1$ and $l,m\in N^+$, we have
\begin{equation*}
| {{n_{m+l}} - {n_{m}}}| \leq {K(\eta) }\sum_{i=m}^{m+l-1}\omega_{\Omega}(\eta^i) ,
\end{equation*}
where $n_k$ denotes ${n_{{\Gamma _{{\eta ^k}}}}}$ for convenience. Since $\omega_{\Omega}\left( r \right)$ is a Dini function, there exists a unit vector $n_0$ satisfying
\begin{equation*}
\mathop {\lim }\limits_{k \to \infty } {{n_{k}} = n_0}.
\end{equation*}

Now, for any $\varepsilon>0$, there exists $k_0\geq 0$ such that $|n_k-n_0|\leq \varepsilon/2$ for any $k\geq k_0$ and $\sup_{\theta\in[\eta,1]}K(\theta)\cdot \omega_{\Omega}(\eta^{k_0})\leq \varepsilon/2$. Then for any $0<r<R\eta^{k_0}$, there exists $k\geq k_0$ such that $R\eta^{k+1}\leq r<R\eta^{k}$. Then for some $\eta\leq \theta<1$,
\begin{equation*}
  |n_{\Gamma_r}-n_0|\leq |n_{\Gamma_r}-n_k+n_k-n_0|\leq K(\theta)\omega_{\Omega}(\eta^{k})+\varepsilon/2\leq \varepsilon.
\end{equation*}
Therefore, $\lim_{r \to 0}n_{\Gamma_r} = n_0$.~\qed
~\\

\begin{remark}\label{r-3}
Without loss of generality, we always assume that $n_0=e_n$ throughout this paper.
\end{remark}

Next, we introduce the following lemma, which concerns the boundary $C^{1,\alpha}$ regularity for solutions with flat boundaries. It was first proved by Krylov \cite{MR688919} and further simplified by Caffarelli (see \cite[Theorem 9.31]{MR1814364} and \cite[Theorem 4.28]{MR787227}). We will use the solutions in this lemma to approximate the solutions in \Cref{t-2} and \Cref{t-H-2}.
\begin{lemma}\label{l-1}
Let $u$ satisfy
\begin{equation}
\left\{\begin{aligned}
&u\in S(\lambda,\Lambda,0)&& ~~\mbox{in}~~B_1^+;\\
&u=0&& ~~\mbox{on}~~T_1.
\end{aligned}\right.
\end{equation}

Then $u$ is $C^{1,\alpha}$ at $0$ and
\begin{equation*}
  |u(x)-u(0)-ax_n|\leq C |x|^{1+\alpha}\|u\|_{L^{\infty }(B_1^+)}, ~~\forall ~x\in B_{1/2}^+
\end{equation*}
with
\begin{equation*}
  |a|\leq C\|u\|_{L^{\infty }(B_1^+)},
\end{equation*}
where $\alpha$ and $C$ depend only on $n, \lambda$ and $\Lambda$.
\end{lemma}

Now, we give the

\noindent\textbf{Proof of \Cref{t-2}.} Let $\omega(r)=\max \left\{\omega_{\Omega}(r),\omega_g(r),\omega_f(r)\right\}$. From the Dini condition, there exists $r_1>0$ such that for any $0<r\leq r_1$,
\begin{equation}\label{e-H-dini-2}
 \omega(r)\leq c_0 ~~\mbox{and}~~ \int_{0}^{r}\frac{\omega(s)}{s}ds\leq c_0,
\end{equation}
where $c_0\leq 1/4$ is a small constant to be specified later and depends only on $n,\lambda,\Lambda$ and $K$. By a proper scaling, we assume that $r_1=1$. Furthermore, we assume that $u(0)=g(0)=0$ and $\nabla g(0)=0$. Otherwise, we may consider $v:=u-g(0)-\nabla g(0)\cdot x$, which satisfies the same equation.

Let $M=\|u\|_{L^{\infty }(\Omega)}+\|f\|_{L^{n}(\Omega)}+[g]_{C^{1,\mathrm{Dini}}(0)}$ and $\Omega _{r}=\Omega \cap B_{r}$. To prove that $u$ is $C^{0,1}$ at $0$, we only need to prove the following:

There exist constants $0<\alpha_{0}, \eta < 1,\bar{C}$ (depending only on $n$, $\lambda$, $\Lambda$), $\hat{C}$ (depending only on $n,\lambda,\Lambda$ and $K$) and a nonnegative sequence $\{a_k\}$ ($k\geq -1$) such that for all $k\geq 0$
\begin{equation}\label{e1.16}
\sup_{\Omega _{\eta^k}}(u-a_k{n_{k}} \cdot x)\leq \hat{C} M \eta ^{k}A_k
\end{equation}
and
\begin{equation}\label{e1.17}
|a_k-a_{k-1}|\leq \bar{C}\hat{C}MA_k,
\end{equation}
where
\begin{equation}\label{e-Ak}
n_k=n_{\Gamma_{\eta^k}},A_0=c_0, A_k=\max(\omega(\eta^{k}),\eta^{\alpha_0} A_{k-1}) (k\geq 1).
\end{equation}

Indeed, from \cref{e-Ak}, we have for any $k\geq 1$,
\begin{equation*}
  A_k\leq \omega(\eta^k)+\eta^{\alpha_0}A_{k-1}.
\end{equation*}
Hence,
\begin{equation*}
  \sum_{i=0}^{k} A_i\leq \sum_{i=1}^{k}\omega(\eta^i)+\eta^{\alpha_0}\sum_{i=0}^{k} A_i+c_0,
\end{equation*}
which indicates
\begin{equation}\label{e.n1}
\begin{aligned}
\sum_{i=0}^{k} A_i\leq \frac{1}{1-\eta^{\alpha_0}}\left(\sum_{i=1}^{\infty}\omega(\eta^i)+c_0\right)\leq C
\end{aligned}
\end{equation}
for some constant $C$ independent of $k$. That is, $\sum_{i=0}^{\infty} A_i$ converges. Thus, $a_k$ converges to some constant $a$.

Then for any $r>0$, there exists $k\geq 0$ such that $\eta^{k+1}<r\leq \eta^{k}$. From \cref{e1.16} and \cref{e1.17}, we have
\begin{equation*}
\sup_{\Omega _{r}}u\leq \sup_{\Omega _{\eta^k}}u\leq \hat{C} M \eta ^{k}A_k+CM\eta^k\leq CM\eta^k\leq CMr,
\end{equation*}
where $C$ depends only on $n,\lambda,\Lambda$ and $K$. In addition,
\begin{equation*}
\inf_{\Omega _{r}}u\geq -C M r
\end{equation*}
can be proved similarly. Therefore,
\begin{equation*}
  \|u\|_{L^{\infty}(\Omega_{r})}\leq CMr.
\end{equation*}
That is, $u$ is $C^{0, 1}$ at $0$.

Now, we prove \cref{e1.16} and \cref{e1.17} by induction. For $k=0$, by setting $a_{-1}=a_0=0$, they hold clearly provided
\begin{equation}\label{e.21}
\hat{C}c_0\geq 1.
\end{equation}
Suppose that they hold for $k$. We need to prove that they hold for $k+1$.

For convenience, we use the following notations. Let $r=\eta ^{k}$, $B_{{\Gamma _r}}^ + =B_r\cap \left\{x|x\cdot{n_{{\Gamma _r}}}>0\right\}$ and ${T_{{\Gamma _r}}} = {B_r} \cap \left\{ {x\left| {x \cdot {n_{{\Gamma _r}}} = 0} \right.} \right\}$ where ${{\Gamma _r}}$ denotes a hyperplane depending only on $r$ and ${n_{{\Gamma _r}}}$ is the unit normal vector of ${{\Gamma _r}}$. We may also denote $n_{\Gamma _{\eta^k}}(n_{\Gamma _r})$ by $n_k$.

Since $\partial\Omega$ satisfies the exterior Reifenberg ${C^{1,\mathrm{Dini}}}$ condition at $0$, there exist a hyperplane ${\Gamma _r}$ and its unit normal vector $n_{\Gamma _r}$ such that $B_r\cap\Omega\subset B_r\cap\{x\cdot {n_{\Gamma _r}}>-r\omega(r)\}$. Let $\tilde B_{{\Gamma _r}}^ +  = B_{{\Gamma _r}}^ + -K(\eta) r\omega \left( r \right) n_{\Gamma _r},$ $\tilde T_{{\Gamma _r}}^ += T_{{\Gamma _r}}^ +  - K(\eta)r\omega \left( r \right) n_{\Gamma _r}$ and $\tilde{\Omega }_{r}=\Omega \cap \tilde B_{{\Gamma _r}}^ + $. Take
\begin{equation}\label{e.K}
c_0K(\eta)<\frac{1}{4}.
\end{equation}
Then $\omega(r) \leq \omega(1)\leq c_0\leq 1/(4K(\eta))$ and $\Omega _{\eta r}\subset\tilde{\Omega }_{r}$.

Let $v$ solve
\begin{equation*}
\left\{\begin{aligned}
 &M^{+}(D^2v,\lambda,\Lambda)=0 &&\mbox{in}~~\tilde{B}^{+}_{{\Gamma _r}}; \\
 &v=0 &&\mbox{on}~~\tilde{T}_{{\Gamma _r}};\\
 &v=\hat{C} M \eta ^{k}A_k &&\mbox{on}~~\partial \tilde{B}^{+}_{{\Gamma _r}}\backslash \tilde{T}_{{\Gamma _r}}.
\end{aligned}
\right.
\end{equation*}
Let $w=u-a_k{n_{k}} \cdot x-v$. Then $w$ satisfies (note that $v\geq 0$ in $\tilde{B}^{+}_{{\Gamma _r}}$)
\begin{equation*}
    \left\{
    \begin{aligned}
      &w\in \underline{S}(\lambda /n,\Lambda , f) &&\mbox{in}~~ \Omega \cap \tilde{B}^{+}_{{\Gamma _r}}; \\
      &w\leq g-a_k {n_{k}} \cdot x &&\mbox{on}~~\partial \Omega \cap \tilde{B}^{+}_{{\Gamma _r}};\\
      &w\leq 0 &&\mbox{on}~~\partial \tilde{B}^{+}_{{\Gamma _r}}\cap \bar{\Omega}.
    \end{aligned}
    \right.
\end{equation*}

In the following arguments, we estimate $v$ and $w$ respectively. By the boundary $C^{1,\alpha}$ estimate for $v$ (see \Cref{l-1}) and the maximum principle, there exist $0<\alpha<1$ (depending only on $n,\lambda$ and $\Lambda$) and $\bar{a}\geq 0$  such that
\begin{equation*}
\begin{aligned}
    \|v-\bar{a}({n_{k}}\cdot x+K(\eta)r\omega \left( r \right){n_{k}})\|_{L^{\infty }(\Omega _{\eta r})}&\leq C_1\frac{(\eta r)^{1+ \alpha}}{r^{1+ \alpha}}\|v\|_{L^{\infty }( \tilde{B}^{+}_{\Gamma_r})}\\
    &\leq C_1\eta ^{\alpha-\alpha_0 }\cdot \hat{C}M\eta ^{(k+1)}\eta^{\alpha_0}A_k\\
     &\leq C_1\eta ^{\alpha-\alpha_0 }\cdot \hat{C}M\eta ^{(k+1)}A_{k+1}
\end{aligned}
\end{equation*}
and
\begin{equation}\label{e.19}
0\leq\bar{a}\leq C_2\hat{C}MA_k,
\end{equation}
where $C_1$ and $C_2$ depend only on $n,\lambda$ and $\Lambda$. Take $\alpha_0=\alpha/2$ and then
\begin{equation}\label{e.20}
\begin{aligned}
\|v-\bar{a} n_{k}\cdot x\|_{L^{\infty }(\Omega _{\eta^{k+1}})}=&\|v-\bar{a} n_{k}\cdot x\|_{L^{\infty }(\Omega _{\eta r})}\\
\leq & C_1\eta ^{\alpha_0 }\cdot \hat{C}M\eta ^{(k+1)}A_{k+1}+\bar{a}K(\eta)r\omega(r)\\
\leq &\left( C_1\eta ^{\alpha_{0} }+\frac{C_2K(\eta)\omega(\eta^{k})}{\eta^{1+\alpha_{0}}}\right)\cdot \hat{C}M\eta ^{(k+1)}A_{k+1}\\
\leq &\left( C_1\eta ^{\alpha_{0} }+\frac{C_2K(\eta)c_0}{\eta^{1+\alpha_{0}}}\right)\cdot \hat{C}M\eta ^{(k+1)}A_{k+1}.
\end{aligned}
\end{equation}

For $w$, by the Alexandrov-Bakel'man-Pucci maximum principle, we have
\begin{equation*}
  \begin{aligned}
\sup_{\Omega_{\eta^{k+1}}}w\leq \sup_{\tilde{\Omega} _{r}}w& \leq\|g\|_{L^{\infty }(\partial \Omega \cap \tilde{B}^{+}_{\Gamma _{ r}})}+\sup_{\partial \Omega \cap \tilde{B}^{+}_{\Gamma _r}}(-a_k n_{k} \cdot x)+C_3r\|f\|_{L^n(\tilde{\Omega}_{r})}\\
    &\leq Mr\omega_g(r) +\sum_{i=0}^{k}|a_i-a_{i-1}|\eta^k\omega(\eta^k)+C_3r\|f\|_{L^{n}(\Omega)}\omega_f(r)\\
    &\leq M \eta^k \omega(\eta^k)+\bar{C}\hat{C}M\sum_{i=0}^{k}A_i\eta^k\omega( \eta^k)+C_3M\eta^k \omega(\eta^k),
  \end{aligned}
\end{equation*}
where $C_3$ depends only on $n,\lambda$ and $\Lambda$.

From \cref{e.n1}, we have
\begin{equation}\label{e.22}
\begin{aligned}
  \sum_{i=0}^{\infty} A_i&\leq \frac{1}{1-\eta^{\alpha_0}}\sum_{i=1}^{\infty}\omega(\eta^i)+c_0\\
&=\frac{1}{1-\eta^{\alpha_0}}\sum_{i=1}^{\infty}\frac{\omega(\eta^i)
\left(\eta^{i-1}-\eta^i\right)}{\eta^{i-1}-\eta^i}+c_0\\
&=\frac{1}{\left(1-\eta^{\alpha_0}\right)\left(1-\eta\right)}\sum_{i=1}^{\infty}
\frac{\omega(\eta^i)\left(\eta^{i-1}-\eta^i\right)}{\eta^{i-1}}+c_0\\
&\leq \frac{1}{\left(1-\eta^{\alpha_0}\right)\left(1-\eta\right)}\int_{0}^{1}
\frac{\omega(r)dr}{r}+c_0\\
&\leq \frac{c_0}{\left(1-\eta^{\alpha_0}\right)\left(1-\eta\right)}+c_0\leq 3c_0,
\end{aligned}
\end{equation}
provided
\begin{equation}\label{e-w}
\left(1-\eta^{\alpha_0}\right)\left(1-\eta\right)\geq 1/2.
\end{equation}
From the definition of $A_k$ again,
\begin{equation*}
  \omega(\eta^{k})\leq A_k\leq \frac{A_{k+1}}{\eta^{\alpha_0}}.
\end{equation*}
Hence,
\begin{equation}\label{e1.22}
  \begin{aligned}
\sup_{\Omega_{k+1}}w &\leq \frac{1}{\eta^{1+\alpha_0}}
M\eta^{k+1}A_{k+1}+\frac{3c_0\bar{C}}{\eta^{1+\alpha_0}}\hat{C}M\eta^{k+1}A_{k+1}+\frac{C_3 }{\eta^{1+\alpha_0}} M\eta^{k+1}A_{k+1}\\
&\leq \frac{C_3+1}{\eta^{1+\alpha_0}}
M\eta^{k+1}A_{k+1}+\frac{3c_0\bar{C}}{\eta^{1+\alpha_0}}\hat{C}M\eta^{k+1}A_{k+1}\\
&\leq\left(\frac{C_3+1}{\hat{C}\eta^{1+\alpha_0}}+\frac{3c_0\bar{C}}{\eta^{1+\alpha_0}}\right)
\hat{C}M\eta^{k+1}A_{k+1}.
\end{aligned}
\end{equation}
Since
\begin{equation*}
  \begin{aligned}
\left| n_{k} - n_{k+1}\right| \leq {K(\eta) }\omega \left( {{\eta ^k}} \right),
\end{aligned}
\end{equation*}
we have for $x\in B_{\eta^{k+1}}$,
\begin{equation}\label{e-H-1.25}
  \begin{aligned}
\left| {{a_{k + 1}}\left(n_k-n_{k+1} \right) \cdot x} \right| &\leq {a_{k + 1}}{K(\eta) }\omega ( \eta ^k){\eta ^{k + 1}}\\
 &\leq \frac{{{a_{k + 1}}{K(\eta) }}}{{\hat CM{\eta ^{{\alpha _0}}}}}\hat C M{\eta ^{\left( {k + 1} \right)}}{A_{k + 1}}\\
 &\leq \frac{3c_0 K(\eta)\bar{C}}{\eta^{\alpha_0}}\hat C M{\eta ^{\left( {k + 1} \right)}}{A_{k + 1}}.
\end{aligned}
\end{equation}

Let $\bar{C}=C_2/{\eta ^{{\alpha _0}}}$. Take $\eta $ small enough such that \cref{e-w} holds and
\begin{equation*}
  C_1\eta ^{\alpha_0 }\leq \frac{1}{6}.
\end{equation*}
Take $c_0$ small enough such that
\begin{equation*}
\frac{C_2K(\eta)c_0}{\eta^{1+\alpha_0}}\leq\frac{1}{6},
\frac{3c_0\bar{C}}{\eta^{1+\alpha_0}}\leq\frac{1}{6}~\mbox{and}~
\frac{3c_0 K(\eta)\bar{C}}{\eta^{\alpha_0}}\leq \frac{1}{3}.
\end{equation*}
Finally, take $\hat{C}$ large enough such that \cref{e.21} holds and
\begin{equation*}
  \frac{C_3+1}{\hat{C} \eta ^{1+\alpha_0}}\leq \frac{1}{6}.
\end{equation*}

Let $a_{k+1}=a_k+\bar{a}$. Then combining \cref{e.20}, \cref{e1.22} and \cref{e-H-1.25}, we have for $x\in \Omega_{\eta^{k+1}}$,
\[\begin{array}{l}
\quad u - {a_{k + 1}}{n_{k + 1}} \cdot x = u - {a_{k + 1}}{n_{k}} \cdot x + {a_{k + 1}}{n_{k}} \cdot x - {a_{k + 1}}{n_{k + 1}} \cdot x\\
 = u - \left( {{a_k} + \bar a} \right){n_{k}} \cdot x + {a_{k + 1}}\left( {{n_{k}} - {n_{k + 1}}} \right) \cdot x\\
 = u - {a_k}{n_{k}} \cdot x - v + v - \bar a{n_k} \cdot x + {a_{k + 1}}\left( {{n_{k}} - {n_{k + 1}}} \right) \cdot x\\
 \le \hat CM{\eta ^{\left( {k + 1} \right)}}{A_{k + 1}}.
\end{array}\]
By induction, the proof is completed. \qed~\\

The proof of the Hopf lemma is similar to that of the boundary Lipschitz regularity. Here, we focus on the curved boundary toward the interior of the domain. We need the following lemma, which can be easily proved by constructing a proper barrier.
\begin{lemma}\label{le-H-1}
Let $u\geq 0$ satisfy
\begin{equation*}
  \left\{\begin{aligned}
    M^{-}(D^2u,\lambda,\Lambda)&=0~~\mbox{in}~~Q_1^+;\\
    u&=0~~\mbox{on}~~T_1;\\
    u&\geq 1~~\mbox{on}~~T_1+e_n.
\end{aligned}\right.
\end{equation*}
Then
\begin{equation*}\label{e-H-l}
u(x)\geq c_1x_n ~~\mbox{in}~~B_{\delta_1}^+,
\end{equation*}
where $\delta_1>0$ and $c_1>0$ depend only on $n,\lambda$ and $\Lambda$.
\end{lemma}

Now, we give the~\\
\noindent\textbf{Proof of \Cref{t-H-2}.} As before, from the Dini condition, there exists $0<r_1<R$ such that for any $0<r\leq r_1$,
\begin{equation}\label{e-H-dini-4}
 \omega(r)\leq c_0 ~~\mbox{and}~~ \int_{0}^{r}\frac{\omega(s)}{s}ds\leq c_0,
\end{equation}
where $c_0\leq 1/4$ is a small constant to be specified later and depends only on $n,\lambda,\Lambda$ and $K$. Moreover, since $n_{\Gamma_r}$ converges as $r\to 0$ (see \Cref{t-1}), we assume that $n_{\Gamma_r} \to e_n$ without loss of generality. Then there exists $0<r_1<R$ such that for any $0<r\leq r_2$,
\begin{equation}\label{e-H-dini-3}
\left| {{n_{\Gamma_r}} - {e_n}} \right| \leq c_0.
\end{equation}
By a proper scaling, we assume that $\min\{r_1,r_2\}=1$ and $u(e_n/2)=1$. Let $\Omega _r^ +  = \Omega  \cap B_{{\Gamma _r}}^ +$. To prove \cref{e-H-main}, we only need to prove the following:

There exist constants $0<\alpha_{0}, \eta < 1$, $\bar{C}$ and $\tilde{a}>0$ (depending only on $n,\lambda$ and $\Lambda$), $\hat{C}$ (depending only on $n,\lambda,\Lambda$ and $K$) and a nonnegative sequence $\{a_k\}$ such that for $k\geq0$,
\begin{equation}\label{e-H1.16}
\inf_{\Omega _{{\eta ^{k + 1}}}^ + }(u-\tilde{a}{n_{k}}\cdot x+a_k{n_{k}}\cdot x)\geq -\hat{C} \eta ^{k}A_k,
\end{equation}
\begin{equation}\label{e-H1.17}
|a_k-a_{k-1}|\leq \bar{C}\hat{C}A_k
\end{equation}
and
\begin{equation}\label{e-H-atilde}
a_k\leq \frac{\tilde{a}}{2},
\end{equation}
where $A_k$ and $n_k$ are defined as before.

Indeed, for any unit vector $l\in R^n$ with $l\cdot {e_n}=\tau>0$, there exists $k_0\geq1$ such that $\omega(\eta^{k})\leq \tau/4$ and $l\cdot n_k\geq \tau/2$ for any $k\geq k_0$. Then $tl\in \Omega$ for any $0<t<1$. Note that $A_k\to 0$ as $k\to \infty$ and then there exits $k_1\geq k_0$ such that
\begin{equation*}
  A_k\cdot \frac{\hat{C}}{ l\cdot{n_{k}}\eta^2}\leq \frac{\tilde{a}}{4}.
\end{equation*}
Take $\delta=\eta^{k_1}$. For $0<t<\delta$, there exists $k\geq k_1$ such that
\begin{equation*}
  \eta^{k+2}\leq t \leq \eta^{k+1}.
\end{equation*}
Then by \cref{e-H1.16},
\begin{equation*}
\begin{aligned}
u(tl)&\geq \tilde{a}({n_{k}}\cdot l)t-a_k({n_{k}}\cdot l)t-\hat{C}\eta^kA_k
&\geq \frac{\tilde{a}({n_{k}}\cdot l)t}{2}-\frac{\hat{C}A_kt}{\eta^2}
&\geq \frac{\tilde{a}({n_{k}}\cdot l)t}{4}
&\geq \frac{\tilde{a}\tau t}{4}.
\end{aligned}
\end{equation*}
That is, \cref{e-H-main} holds.

Now, we prove \crefrange{e-H1.16}{e-H-atilde} by induction. Let $\tilde{Q}=Q^{+}_{1/2}+c_0{e_n},\tilde{T}=\partial \tilde{Q}\cap \left\{x_n=1/4+c_0\right\}$. That is, $\tilde{T}$ is the top boundary of $\tilde{Q}$. Note that $c_0\leq 1/4$. Thus, $\tilde{Q}\subset \Omega\cap B^+_{1}$.

By the Harnack inequality,
\begin{equation*}
  \inf_{\tilde{T}} u\geq \tilde{c} u({e_n}/2)=\tilde{c},
\end{equation*}
where $\tilde{c}$ depends only on $n,\lambda$ and $\Lambda$. Let $\tilde{u}$ solve
\begin{equation*}
  \left\{\begin{aligned}
&M^{-}(D^2\tilde{u},\lambda,\Lambda)= 0 &&\mbox{in}~~\tilde{Q}; \\
&\tilde{u}=\tilde{c}&&\mbox{on}~~\tilde{T};\\
&\tilde{u}=0 &&\mbox{on}~~\partial \tilde{Q}\backslash \tilde{T}.
  \end{aligned}\right.
\end{equation*}
From \Cref{le-H-1}, there exist $\delta_1>0$ and $0<c_2<1/2$ such that
\begin{equation*}
u(x)\geq \tilde{u}(x) \geq c_2(x_n-c_0) ~~\mbox{in}~~B^+_{\delta_1}+c_0{e_n}.
\end{equation*}
Note that $u\geq 0$ and hence,
\begin{equation*}
  u(x)\geq c_2(x_n-c_0) ~~\mbox{in}~~\Omega\cap B_{\delta_1}.
\end{equation*}
Therefore, by noting \cref{e-H-dini-3},
\begin{equation*}
  u(x)\geq c_2x\cdot n_0-2c_2c_0\geq c_2x\cdot n_0-c_0 ~~\mbox{in}~~\Omega\cap B_{\delta_1}.
\end{equation*}

Set $\tilde{a}=c_2$ and $a_{-1}=a_0=0$. Take
\begin{equation}\label{e-H-1}
\eta\leq \delta_1.
\end{equation}
Then \crefrange{e-H1.16}{e-H-atilde} hold for $k=0$. Suppose that they hold for $k$. We need to prove that they hold for $k+1$.

Let $r=\eta ^{k+1},\tilde{B}^+=B^{+}_{\Gamma _r}-K(\eta)r\omega(r)n_{\Gamma_r}$ and $\tilde{T}=T_{\Gamma _r}-K(\eta)r\omega(r)n_{\Gamma_r}$. Then $B^{+}_{\Gamma _{\eta r}}\subset \tilde{B}^+$. Let $v$ solve
\begin{equation*}
\left\{\begin{aligned}
 &M^{-}(D^2v,\lambda,\Lambda)=0 &&\mbox{in}~~\tilde{B}^+; \\
 &v=0 &&\mbox{on}~~\tilde{T};\\
 &v=-\hat{C}\eta ^{k}A_k &&\mbox{on}~~\partial \tilde{B}^{+}\backslash T_{\Gamma _r}.
\end{aligned}
\right.
\end{equation*}
Let $w=u-\tilde{a}{n_{k}}\cdot x+a_k{n_{k}}\cdot x-v$. Then $w$ satisfies (note that $u\geq 0$ and $v\leq 0$)
\begin{equation*}
    \left\{
    \begin{aligned}
      &M^{-}(D^2w,\lambda/n,\Lambda)\leq 0 &&\mbox{in}~~ \Omega\cap \tilde{B}^{+}; \\
      &w\geq -\tilde{a}{n_{k}}\cdot x+a_k{n_{k}}\cdot x&&\mbox{on}~~\partial \Omega \cap \tilde{B}^{+};\\
      &w\geq 0 &&\mbox{on}~~\partial  \tilde{B}^{+}\cap \bar{\Omega}.
    \end{aligned}
    \right.
\end{equation*}

In the following arguments, we estimate $v$ and $w$ respectively. By the boundary $C^{1,\alpha}$ estimate for $v$ (see \Cref{l-1}) and the maximum principle, there exist $0<\alpha<1$ (depending only on $n,\lambda$ and $\Lambda$) and $\bar{a}\geq 0$  such that (note that $A_k\leq A_{k+1}/\eta^{\alpha_0}$)
\begin{equation*}
\begin{aligned}
    \|v&+\bar{a}{n_{k+1}}\cdot (x-K(\eta)r\omega(r)n_{k+1})\|_{L^{\infty }(B^+ _{\Gamma _{\eta r}})}\leq C_1 \eta ^{1+ \alpha}\|v\|_{L^{\infty }(\tilde{B}^{+})}\\
    &\leq C_1\eta ^{\alpha}\cdot \hat{C}\eta ^{(k+1)}A_k
    \leq C_1\eta ^{\alpha-\alpha_0 }\cdot \hat{C}\eta ^{(k+1)}A_{k+1}
\end{aligned}
\end{equation*}
and
\begin{equation}\label{e.H-19}
\bar{a}\leq C_2\hat{C}A_k/\eta,
\end{equation}
where $C_1$ and $C_2$ depend only on $n,\lambda$ and $\Lambda$. Take $\alpha_0=\alpha/2$. Then
\begin{equation}\label{e.H-20}
\begin{aligned}
\|v+\bar{a}{n_{k+1}}\cdot x\|_{L^{\infty }(\Omega^+_{k+2})}&\leq\|v+\bar{a}{n_{k+1}}\cdot x\|_{L^{\infty }(B^+_{\Gamma _{\eta r}})}\\
&\leq \left(C_1\eta ^{\alpha_0 }+\frac{C_2A_kK(\eta)}{\eta}\right)\cdot \hat{C}\eta ^{(k+1)}A_{k+1}\\
&\leq \left(C_1\eta ^{\alpha_0 }+\frac{C_2c_0K(\eta)}{\eta}\right)\cdot \hat{C}\eta ^{(k+1)}A_{k+1}
\end{aligned}
\end{equation}

For $w$, by the maximum principle, we have
\begin{equation*}
  \begin{aligned}
\inf_{\Omega^+_{k+2}}w&\geq \inf_{\Omega\cap\tilde{B}^+}w \geq \inf_{\partial \Omega \cap \tilde{B}^{+}}\left(-\tilde{a}{n_{k}}\cdot x+a_k n_{k}\cdot x\right)\\
&\geq \inf_{\partial \Omega \cap \tilde{B}^{+}}\left(-\tilde{a}|n_{k}-n_{k+1}|\cdot|x|-\tilde{a}(n_{k+1}\cdot x)\right) \\
&\geq -\tilde{a}K(\eta)\eta^{k + 1}\omega(\eta^k)-\tilde{a}{\eta ^{k + 1}}\omega(\eta ^{k + 1}).
  \end{aligned}
\end{equation*}
As before, from the definition of $A_k$ (see \cref{e-Ak}),
\begin{equation*}
  \omega(\eta^{k+1})\leq A_{k+1}~\mbox{and}~\omega(\eta^{k})\leq A_{k}\leq \frac{A_{k+1}}{\eta^{\alpha_0}}.
\end{equation*}
Hence,
\begin{equation}\label{e-H-1.21}
  \begin{aligned}
\inf_{\Omega^+_{k+2}}w &\geq -\frac{K(\eta)+1}{\hat{C}\eta^{\alpha_0}}\cdot \hat{C}\eta^{k+1}A_{k+1}.
\end{aligned}
\end{equation}

Since $\left| {{n_{k}} - {n_{k + 1}}} \right| \le {K(\eta) }\omega (\eta ^k)$, for $x\in B_{\eta^{k+2}}$,
\begin{equation}\label{e-H-1.22}
  \begin{aligned}
\left| {\tilde a\left( {{n_{k}} - {n_{k + 1}}} \right) \cdot x} \right|
  \leq \tilde a{K(\eta) }\omega \left( {{\eta ^k}} \right)\eta^{k+2}
  \leq \frac{K(\eta) }{\hat C}\hat C{\eta ^{\left( {k + 1} \right)}}{A_{k + 1}}.
\end{aligned}
\end{equation}

Take $\eta $ small enough such that \cref{e-H-1} holds,
\begin{equation*}
  C_1\eta ^{\alpha_0 }\leq 1/6
\end{equation*}
and
\begin{equation}\label{e-H-2}
\left(1-\eta^{\alpha_0}\right)\left(1-\eta\right)\geq 1/2.
\end{equation}

Let $\bar{C}=C_2/{{\eta ^{1 + {\alpha _0}}}}$. Take $\hat{C}$ large enough such that
\begin{equation*}
\frac{K(\eta)+1}{\hat C\eta ^{\alpha _0}} \leq \frac{1}{3}.
\end{equation*}
As before, by noting that \cref{e.22} and \cref{e-H-2}, we have
\begin{equation*}
  \sum_{i=0}^{\infty} A_i\leq 3c_0.
\end{equation*}
Finally, take $c_0$ small enough such that
\begin{equation*}
3c_0 \bar{C}\hat{C}\leq\frac{\tilde{a}}{2}~\mbox{and}~\frac{C_2c_0K(\eta)}{\eta} \leq \frac{1}{6}.
\end{equation*}

Let $a_{k+1}=a_k+\bar{a}$. Then combining \crefrange{e.H-20}{e-H-1.22}, we have for $x\in \Omega^+_{k+2}$
\begin{equation*}
  \begin{aligned}
u &-\tilde a{n_{k+1}} \cdot x + {a_{k + 1}}{n_{k+1}} \cdot x\\
&=u-\tilde a n_{k} \cdot x + a_{k} n_{k} \cdot x-\tilde a(n_{k+1}-n_{k}) \cdot x
+ a_k(n_{k+1}-n_{k}) \cdot x+ \bar{a}n_{k+1} \cdot x\\
&= u - \tilde a{n_{k}} \cdot x + {a_k}{n_{k}} \cdot x - v + \left( {v + \bar{a}n_{k+1} \cdot x} \right) + (a_k-\tilde a)\left( {{n_{k}} - {n_{k+1}}} \right) \cdot x\\
&= \omega  + \left( {v + \bar{a}n_{k+1} \cdot x} \right)+ (a_k-\tilde a)\left( {{n_{k}} - {n_{k+1}}} \right) \cdot x\\
&\geq  - \hat C{\eta ^{\left( {k + 1} \right)}}{A_{k + 1}}
  \end{aligned}
\end{equation*}
and
\begin{equation*}
  a_{k+1}\leq \sum_{i=0}^{k}|a_i-a_{i-1}|\leq\bar{C}\hat C\sum\limits_{i = 0}^\infty  {{A_i}}\leq 3c_0 \bar{C}\hat{C}\leq \frac{\tilde{a}}{2}.
\end{equation*}
By induction, the proof is completed. \qed~\\

\bibliographystyle{model4-names}
\bibliography{L_H}



\end{document}